\documentclass{amsart}

\usepackage{amsbsy}
\usepackage{graphicx}
\graphicspath{{converted_graphics/}}
\begin{document}

\begin{center}
\textbf{\LARGE \quad\quad          Asymptotic Admissibility of Priors\newline and Elliptic Differential Equations}
\end{center}

\begin{center}
\Large J.A.Hartigan ,Yale University
\end{center}

\textbf{Abstract }
\textit{
We evaluate priors by the second order asymptotic behaviour of the corresponding estimators.
Under certain regularity conditions, the risk differences between efficient 
estimators of parameters taking values in a domain $D$, an open connected subset of 
$R^d$, are asymptotically expressed as elliptic differential forms 
depending on the asymptotic covariance matrix $V$. Each efficient 
estimator has the same asymptotic risk as a ``local Bayes'' estimate 
corresponding to a prior density $p$. The asymptotic decision theory of the 
estimators identifies the smooth prior densities as admissible or 
inadmissible, according to the existence of solutions to certain elliptic 
differential equations. The prior $p$ is admissible if the quantity $pV$ is sufficiently small 
near the boundary of $D$. We exhibit the unique admissible invariant prior for $V=I,D=R^d-\{0\}$.
 A detailed example is given for a normal mixture model.}
\newline \newline
\textbf{1 Introduction}
\newline \newline
A parameter $x$ takes values in a \textit{domain} $D$, an open connected subset of $R^{d}$. \newline
I use the partial differential equation symbol $x$ rather than the usual statistical 
symbol $\theta $ because the evaluation of asymptotic risk 
reduces to existence problems in the theory of partial differential 
equations. The parameter indexes a probability density $p(y_n \vert 
x)$ with respect to some measure $\mu _n $, say, for data $y_n $. 
\newline\newline
We use the Kullback-Leibler loss function 
\begin{equation}
\hat {x},x\in D:L_n (\hat {x},x)=\int {\log } [p(y_n \vert x)/p(y_n \vert 
\hat {x})]p(y_n \vert x)d\mu _n (y_n )
\end{equation}
to define the risk of the estimator $\hat {x}_n$, a function of $y_n $ taking 
values in $ D $:
\begin{equation}
R_n (\hat {x}_n ,x)=\int {L_n (\hat {x}_n (y_n ),x)p(y_n \vert x)d\mu _n } 
.
\end{equation}
For a prior density p, the posterior Bayes estimator $\hat {x}(y_n 
,p)$ minimizes the posterior Kullback-Leibler risk
\begin{equation}
R(\hat {x}\vert y_n )=\int {L(\hat {x},x)p(y_n \vert x)p(x)dx/} \int {p(y} 
_n \vert x)p(x)dx.
\end{equation}
Define 
\begin{equation}
V_n (x)=-1/\int \textstyle{\partial \over {\partial x}}\textstyle{\partial \over {\partial x'}}\log [p(y_n \vert x)]p(y_n \vert 
x)d\mu _n (y_n ). 
\end{equation}
We will assume $nV_n (x)\to V(x), $ the \textit{asymptotic covariance matrix. }
\newpage
\noindent Following Brown[Br79], and letting $U$ denote the prior uniform over D, 
we consider estimators of form $\hat {x}(y_n ,U)+V_n b(\hat 
{x}(y_n ,U))$ for fixed \textit{decision functions} $b:D\to D.$ 
Under smoothness conditions [Ha98] requiring smooth variation of the data density and the prior with 
$x$, the asymptotic risks for different decision functions  differ only by 
terms of order $n^{-2}$; therefore we define the \textit{asymptotic risk} for decision function $b$, 
relative to the decision function $b=0 $ corresponding to the uniform prior $U$, by the assumed limit
\begin{equation}
R(b,x)=\mathop {\lim }\limits_{n\to \infty } n^2[R_n (b,x)-R_n (U,x)]=
\sum\nolimits_{i,j} {\{\partial _i } (V_{ij} b_j)+\textstyle{1 \over 2}b_i 
b_j V_{ij} \}
\end{equation}
where $\partial _i $ denotes the partial derivative $\textstyle{\partial 
\over {\partial x_i }}$.
\newline\newline
For a prior with density $p$, the posterior bayes estimate corresponds 
asymptotically to the decision function $b_i^p =\partial _i \log p, $ and then 
the risk may be expressed in elliptic operator form 
\begin{equation}
R(b^p)=2\sum\nolimits_{ij} {\partial _i } (V_{ij} \partial _j \sqrt p )/\sqrt p . 
\end{equation}
\newline
It turns out that, under certain conditions of smoothness and boundedness 
for $b$ and $V$, there is a\textit{ risk matching} prior density p for which $R(b^p)=R(b). $ Thus the 
behaviour of asymptotic risk for all smooth decisions is captured in the 
theory of elliptic differential equations, equations whose relevance to 
decision theory were first indicated in Stein[St56], but which were 
extensively elucidated for the normal location problem in Brown[Br71].
See also Strawderman and Cohen[SC71]. The 
asymptotic behavior of Bayes estimators near maximum likelihood estimators 
have been studied for loss functions of form $L_n (\hat {x},x)=W(\hat 
{x}-x)$ by Levit in [Le82], [Le83], [Le85]; in particular, he shows that the Bayes 
estimators form a complete class under certain regularity conditions.
\newline \newline
\noindent \textbf{2 Risk matching priors}
\newline \newline
For each decision function $b$ we will find a risk matching prior p for 
which $R(b^p)=R(b)$. Then we need only consider decision functions of form
$b_i^p =\partial _i \log p$ \newline and risks of form $R(b^p)=2\sum\nolimits_{ij} 
{\partial _i } (V_{ij} \partial _j \sqrt p )/\sqrt p $ in the asymptotic 
decision theory. This result will be proved under boundedness and smoothness 
assumptions using some standard tools from Pinsky[Pi95].
\newline\newline
For the domain $D$ with closure $\bar {D}$ , a function f is uniformly 
Holder continuous with exponent $\alpha , 0<\alpha \le 1,$ in $\bar D $ if
\begin{equation}
\vert \vert f\vert \vert _{0,\alpha ,\bar D} =\mathop {\sup 
}\limits_{ x,y\in \bar D , x\ne y } \frac{\vert f(x)-f(y)\vert }{\vert x-y\vert ^\alpha }<\infty 
.
\end{equation}
The Holder spaces $C^{k,\alpha }(\bar {D})$ consist of functions whose $k$-th 
order partial derivatives are uniformly Holder continuous with exponent $\alpha$ in $\bar D$.
 Say $D' \subset \subset D$ if $D'$ is bounded and properly included in $D$. The Holder spaces $C^{k,\alpha }(D)$ consist of functions that lie in $C^{k,\alpha 
}(\bar {D}')$ for each $D' \subset \subset D$.

The domain $D$ has a $C^{k,\alpha }$ boundary $\partial D$ if for each point $x_0 \in 
\partial D$, there is a ball B centered at $x_0 $ and a 1-1 mapping $\psi 
:B\to A\subset R^d$, such that $\psi (B\cap D)\subset \{x\in R^d,x_n 
>0\},\psi (B\cap \partial D)\subset \{x\in R^d,x_n =0\},\psi \in 
C^{k,\alpha }(B),\psi ^{-1}\in C^{k,\alpha }(A).$ 
\newline
\begin{eqnarray}
 \mbox{ Condition }\bar{A}&:& \mbox{ D bounded} , \partial D\in C^{2,\alpha },V_{ij} \in 
C^{1,\alpha }(\bar {D}),V \mbox{ positive definite in }  \bar {D}.\nonumber\\
 \mbox{Condition } A&:& V_{ij} \in C^{1,\alpha }(D), V \mbox{ positive definite in }D.\nonumber
\end{eqnarray}  
 
\noindent We follow a standard approach which first proves results under the strong condition $\bar A$,
and then extends the results to the weak condition $A$ by approximating $D$ by an increasing sequence of bounded subdomains $D_n$.
\newline 

\noindent\textbf{Theorem 1}: { Suppose $ b_i = \partial_i b \in C^{1,\alpha }(\bar D)$  and assumption $\bar A$ holds. 
\newline 
Then there exists a prior $p,\sqrt p  \in C^{2,\alpha} (\bar D),p>0$ in $D,$  such that}
\begin{equation} 
 R(b)=R(b^p)=2\sum\nolimits_{ij} \partial _i (V_{ij} \partial _j \sqrt p)/\sqrt p. 
\end{equation}

\noindent\textbf{Proof}: 

\noindent From [Pi95, theorem 5.5], for some eigenvalue $ \lambda $,  there exists \newline
$ u\in C^{2,\alpha }(\bar {D}),u>0$ in $D,u=0$ in   $\partial D,$  satisfying
 \begin{equation}
 2\sum\nolimits_{ij} {\partial _i } (V_{ij} \partial _j u)-R(b)u=\lambda u. 
\end{equation}
Since
 \begin{eqnarray} 
\int_D {\lambda u^2} &=&\int_D \{ 2\sum\nolimits_{ij} {\partial _i }(V_{ij} \partial _j u)-R(b)u\}u \label{line1}\\ 
&=&\int_D {-\textstyle{1 \over 2}\sum\nolimits_{ij} ( } ub_i -2\partial _i u)(ub_j -2\partial _j u)V_{ij}\le 0,\label{line2}
\end{eqnarray} 
it follows that $ \lambda \le 0. $
If $ \lambda =0, $ the corresponding eigenvector $ u$  provides the $b$\textit{-matching} prior 
$ p=u^2$ with $ b_i =\partial _i \log p.$ 
If $ \lambda <0, $ from [Pi95, theorem 6.5], for each $ \phi \in C^{2,\alpha 
}(\bar {D}),\phi >0$  there exists a unique $b$\textit{-matching} solution \newline
 $ \sqrt p= u\in C^{2,\alpha }(\bar {D}),u>0$ in $D, u=\phi \mbox{ on } \partial D$,  to the 
equation $ R(b)=R(b^p).$ 
\newline\newline
\noindent \textbf{Theorem 2}:  { Suppose $ b_i = \partial_i b \in C^{1,\alpha }( D)$  and assumption A holds. 
\newline 
Then there exists a prior $p,\sqrt p  \in C^{2,\alpha} ( D),p>0$ in $D,$  such that}
\begin{equation} 
 R(b)=R(b^p)=2\sum\nolimits_{ij} \partial _i (V_{ij} \partial _j \sqrt p)/\sqrt p. 
\end{equation}

\noindent\textbf {Proof:} 

\noindent Select an increasing sequence of 
bounded domains $ \bar D_n \subset \subset D_{n+1} ,\cup D_n =D.$  Within each 
domain, assumption $ \bar {A}$  holds, so there exists a sequence of solutions \newline
$ \sqrt {p_n } =u_n \in C^{2,\alpha }(\bar {D}_n ),u_n >0$ in $D_n ,$  such 
that $ R(b)=R(b^{p_n }),x\in D_n .$ 
\newline\newline
For some $ x_0 \in D_1$, without loss of generality  set $ u_n (x_0 )=1,$ all n. Note that any solution 
$ u_N ,N>n$  is also a solution to $ R(b)=R(b^{p_n }),x\in D_n .$ 
A Harnack inequality[Pi95, p 124] implies that  for $N>n,c_n \le u_N 
(x)\le C_n ,x\in D_{n} $ for some constants $c_n,C_n$.
\newline\newline
The Schauder interior estimate[Pi95,p86]
  implies, for some constant $C_n$, for all $N>n$,
\begin{equation}
\vert \vert u_N \vert \vert _{2,\alpha ,D_n } =\mathop {\sup 
}\limits_{\begin{array}{l}
 x,y\in D_n \\ 
 x\ne y \\ 
 \end{array}} \sum\nolimits_{i,j} {\frac{\vert \partial _i \partial _j (u_N 
(x)-u_N (y))\vert }{\vert x-y\vert ^\alpha }} \le C_n  .
\end{equation}
Thus, for each n, the sequence $ u_N $ is precompact in the $ \vert \vert u_N 
\vert \vert _{2,0,D_n } $ norm. By diagonalization, there exists a 
subsequence, say $ u_N$,  that converges to $ u$  in $ C^{2,0}(D)$ 
for which 
\begin{equation}
 R(b)=R(b^p)=2\sum\nolimits_{ij} {\partial _i } (V_{ij} \partial 
_j u)/u,x\in D.  
\end{equation}
Finally, we need to show that $ u\in C^{2,\alpha }(D).$ 
First note that $ u\in C^{2,\alpha }(D_n)$, since
\begin{equation}
  \vert \vert u_N \vert \vert _{2,\alpha ,D_n } \le c_n 
 \mbox{ for } N>n\Rightarrow \vert \vert u\vert \vert _{2,\alpha ,D_n } \le c_{n} .
\end{equation}

\noindent For any $ D'\subset \subset D$ , the compact set 
$ \bar {D}'$  is covered by $ \cup D_n =D$ , and therefore by a finite subcovering, and so by a 
particular $D_n $ . Thus $ \vert \vert u\vert \vert _{2,\alpha ,D'} \le c_n (D')$ for every 
$ D' \subset \subset D$  which implies $ u\in C^{2,\alpha }(D).$ 
\newline \newline
\noindent\textbf{3 Explicit matching priors using the Feynman-Kac integral }
\newline \newline
When $ \lambda =0$ in equation (9), the decision function  $ b $ is the posterior bayes 
decision corresponding to the prior $p: b_i =\partial _i \log p$ . We could 
determine the ratio 
$ p(x)/p(y)$  by the integral$ \int\limits_\rho {\sum\nolimits_i {b_i } } dx_i $  
over \textit{any} path $\rho$ connecting $x$ and $y$.
\newline\newline
When $\lambda <0$, so that $b$ is no longer a gradient, it is plausible
to attempt to find an approximating $p$ by averaging these integrals over \textit{all} 
paths between $x$ and $y$. 
\newline\newline
Consider the stochastic differential equation with initial condition $ X_i 
(0)=x:$ 
\begin{equation}
 dX_i (t)=\sum\nolimits_{ij} V_{ij}^{1/2} dW_j (t)+\textstyle{1 \over 2}(b_i 
-\sum\nolimits_j {\partial _j } V_{ij} )dt  
\end{equation}

\noindent We propose the Feynman Kac integral formula to specify a risk matching 
prior:
\begin{equation}
u(x)=E[\exp (-\textstyle{1 \over 2}\int_0^T {\sum\nolimits_i {b_i } } 
{\circ}dX_i )\sqrt{p(X(T))}]
\end{equation}

\noindent where $ \int_0^T {\sum\nolimits_i {b_i } } \circ dX_i $  is the Stratonovitch 
stochastic integral, and T is the time to reach the boundary of D. I suspect 
that the condition $ \lambda <0$  is sufficient for the existence 
of the stochastic process and the integral when assumption $ \bar {A}$ holds . 
When the formula is valid, we see that $ u(x)$ is determined as a weighted 
combination of its boundary values, with the weight at each boundary point 
determined by the path integral $ \exp (-\textstyle{1 \over 
2}\int\limits_\rho {\sum\nolimits_i {b_i } } dx_i ) $  over the various 
paths that reach that particular boundary point. Many different priors risk matching $b$ are available
corresponding to different smooth assignments to the boundary values.

\newpage
\noindent\textbf{4 Decision theory for asymptotic risks}
\newline \newline
We now apply decision theoretic classifications to the asymptotic risk 
formula. Our conclusions about the prior $ p$  will depend only on the domain $ D$  
and the asymptotic variance $ V$ . From now on we will drop the term asymptotic. 
We consider a particular variance $ V $ and the set of risks, real valued 
functions on the domain $D$, corresponding to priors satisfying
\begin{equation}
\mbox{Assumption }B:p \in C^{2,\alpha }(D),V_{ij} \in C^{1,\alpha 
}(D);p,V>0 \mbox{ in } D.
\end{equation}

\noindent The risks will be written $ R(p)=R(b^p)$.

\noindent The posterior bayes decision $ b^p$  is \textit{locally Bayes}: for any alternative decision

$ b^p+v$  where $ v\in C^{1,\alpha }(D)$ and $v=0$  in  $D-D'$, some $D' \subset \subset
D,$ integration by parts shows
\begin{equation}
\int_D {[R(} b^p+v)-R(b^p)]p\ge 0.
\end{equation}
\textbf{Theorem 3: }Under assumption B, with D bounded,
the following conditions are equivalent:

$ b^p$  is \textit{Bayes: }there exists no $ p^{\ast} \ne p$ with $ \int ( R(p^{\ast} )-R(p))p\le 0.$  

$ b^p$  is \textit{Admissible: } there exists no $ p^\ast \ne p$  with $ R(p^\ast )\le R(p),x\in D.$ 

$b^p $  is \textit{Unique Risk:} there exists no $ p^\ast \ne p$  with $ R(p^\ast )=R(p),x\in D.$ 

$ b^p$  is \textit{Brown:} no non-trivial positive $ h $ solves $ \sum\nolimits_{ij} {\partial _i } (pV_{ij} 
\partial _j h)=0.$ 

\noindent \textbf{Proof:}

\noindent \textit{Bayes }implies \textit{Admissible} because $ u\ast \ne u$  violating\textit{ Admissible} also
violates\textit{ Bayes.}
\textit{Admissible }implies \textit{Unique Risk} because $ u\ast \ne u$  violating\textit{ Unique Risk} also
violates\textit{ Admissible.}
\textit{Brown[Br71,1.3.9]} and \textit{Unique Risk} are equivalent, because $ R(p\sqrt h )=R(p)$  if and only if
$
\sum\nolimits_{ij} {\partial _i } (pV_{ij} \partial _j h)=0.
$
\newline\newline
It only remains to show that failure of \textit{Bayes}
implies failure of \textit{Unique Risk}.

\noindent Without loss of generality, assume $p=U$ is uniform,
so that $R(p)=0.$ If $ p$  is not \textit{Bayes, }there exists $ p^\ast \ne p$ with 
\begin{equation}
0 \ge \int  R(p^\ast ) 
=\textstyle{1 \over 2}\int_D \sum\nolimits_{ij} V_{ij}  b_i^\ast b_j^\ast  
+\int_D \sum\nolimits_{ij} \partial_i (V_{ij} b_j^\ast)
\end{equation}
\noindent  Since $ p^\ast \ne p$ and the 
middle integral is positive, then 
\begin{equation}
\int_D \sum\nolimits_{ij} \partial_i (V_{ij} b_j^\ast)=C<0.
\end{equation}
Let $|\partial D|$ be the lebesgue measure of the boundary $\partial D$ when $D$ is smooth and bounded,
let $\tau$ denote the outward pointing normals on the boundary, and note that 
\begin{equation}
\int_{\partial D} \sum_{ij}\{\tau_i b_i^\ast V_{ij}\}=
\int_D \sum\nolimits_{ij} \partial_i (V_{ij} b_j^\ast)=C<0.
\end{equation}
Applying Theorem 6.31 from [GT97], for $D_n \subset \subset D$,
 since $C<0,\sum_{ij}\{\tau_i \tau_j V_{ij}\}>0$,\newline there exists a solution $p_n =u_n^2 \in C^{2,\alpha}(\bar D_n)$ 
to the oblique derivative problem
 \begin{equation}
R(p_n)=0 \mbox{ in } D_n, u_n=2*C\sum_{ij}\{\tau_i \partial_i u V_{ij}\}/|\partial D_n| \mbox{ in } \partial D_n
\end{equation}
so that $\int_D \sum\nolimits_{ij} \partial_i (V_{ij} \partial_i p_n/p_n )=C<0.$

\noindent Repeating the compactness argument of theorem 2 on $D_n \subset \subset D_{n+1}, \cup D_n = D$, there exists $p_0\in C^{2,\alpha}(D)$ with
\begin{equation}
R(p_0)= 0,  \int_D \sum\nolimits_{ij} \partial_i (V_{ij} \partial_i p_0/p_0 )=C<0.
\end{equation}
The first condition states that $p_0$ and $p=U$ have the same risk, and the second condition guarantees that $p_0 \ne p=U$,
so that the unique risk condition fails, as required. 
\newline\newline

\noindent\textbf{5 When is $pV$  Bayes on} \textbf{R}$ ^{d}$ \textbf{?}
\newline \newline
Brown's condition shows that the locally Bayes estimate $b^p$ is Bayes or not depending only on the product $pV$. For example, $b^p$ is Bayes with $V$, if and only if $b^U$ is Bayes with $pV$. We will therefore rephrase the admissibility question in terms of the product $pV$: the prior-scaled covariance matrix $pV$ is Bayes on $D$ if and only if there is no non-trivial solution to Brown's equation. 
\newline \newline
\textbf{Theorem 4}. Let $D = R ^{d}$  , $pV \in C^{1,\alpha}(D),  r=\vert x\vert ,x=rs$   where $ s$  ranges over the surface $S$ of the unit 
sphere. Define $W(R,s)=[\int_1^R \textstyle { 1 \over p}V^{-1} r^{1-d}dr]^{-1}.$
Suppose that, uniformly over $s \in S$,
\begin{equation}
\lim_{R \to \infty}W(R,s)= W(s),\lim_{R \to \infty}W(s) W^{-1}(R,s)W(s)= W(s).
\end{equation}

Then $pV$ is Bayes on $R^d$  only if 
\begin{equation}
   \int_{s\in S}\sum_{ij}s_i s_j W_{ij}(s)ds=0. 
\end{equation} 
\newline
\noindent\textbf{Proof:}
\newline
Let $ Q=\{Q_i ,1\le i\le d\} \in C^{1,\alpha}(D)$ be a \textit{test} function with
 relative risk 
\begin{equation}
 R(b^p+Q)-R(b)= \sum_{ij} \{\partial _i (pV_{ij} Q_j )+\textstyle{1 \over 2}Q_i 
Q_j pV_{ij}\}.
\end{equation}
From theorem 2 for every $Q$ there exists  a prior $q$ with $ R(b^p+Q)=R(b^q)$. Thus $pV$ is Bayes if and only if

\begin{equation}
\int_D \{R(b^p+Q)-R(b^p)\}p=\int_D\sum_{ij}\{ \partial _i (pV_{ij} Q_j )+\textstyle{1 \over 2}Q_i Q_j pV_{ij} \}\ge 0
\end{equation}
\newline
for every \textit{test} $ Q$ where the integral is defined. 
Equivalently, with the test $(pV)^{-1}Q$,
\begin{equation}
t(Q)=\int_D \{R(b^p+(pV)^{-1}Q)-R(b^p)\}p=\int_D \{ \sum_i \partial_i Q_i + \textstyle{1 \over 2}\sum_{ij} Q_i Q_j 
V_{ij}^{-1}/p \}\ge 0
\end{equation}
for every test $Q$ where the integral is defined.
\newline\newline
The possible negative term$ \int_D \sum_i\partial _i Q_i $ is determined by values 
in the neighbourhood of infinity, so  $pV$ being Bayes is determined by the 
behaviour of $pV$ near the infinite boundary. In particular if two functions $pV$ are identical outside a compact subset of $D$,
they have the same admissibility classification. 
\newline We therefore consider a test function that 
is zero inside the unit sphere:\newline
 \begin{equation}Q(rs)=g(r)r^{1-d}q(s),q(s)=-Ws \end{equation}
where $g$ is twice differentiable, $g(r)=0$ for $0 \le r \le 1,0< g(r)\le 1 $ for $ 1<r<2,g(r)=1$ for $r \ge 2$. 

\begin{equation}
\mbox { Let }t(Q,R)=\int_{|x|<R} \{ \sum_i \partial_i Q_i + \textstyle{1 \over 2}\sum_{ij} Q_i Q_j 
V_{ij}^{-1}/p \}r^{d-1}drds=\int_{s \in S} I(R,s) ds.
\end{equation}

Consider the contribution $I(s,R)$ to the test integral for a particular $s$:
 \begin{eqnarray}
 I(s,R)&=&\sum_i s_i Q_i(Rs) R^{d-1} + \int_0^R \{\textstyle{1 \over 2}\sum_{ij} Q_i Q_j 
V_{ij}^{-1}/p \}r^{d-1}dr \label{line1}\\
&=&\sum_i s_i q_i(s)+  {1 \over 2}\sum_{ij} q_i q_j \int_0^R 
g(r)^2 r^{1-d}V_{ij}^{-1}/p dr \label{line2}\\
&\le&-\sum_{ij} s_i s_j W_{ij}+  {1 \over 2}\sum_{ijkl} s_is_j W_{ik} W_{jl} W_{kl}^{-1}(R,s) \label{line3}\\
&\to&-{1 \over 2}\sum_{ij} s_i s_j W_{ij} \mbox { uniformly in } s \mbox{ as } R  \to \infty \label{line4}
\end{eqnarray}
Thus
\begin{equation}
t(Q)=\lim_{R\to \infty}\int_{s \in S} I(s,R) ds <0
\end{equation}
unless
\begin{equation}
   \int_{s\in S}\sum_{ij}s_i s_j W_{ij}(s)ds=0. 
\end{equation}
which shows that the condition in the theorem is necessary for $p$ to be Bayes.
\newline\newline

Failure of the condition in the theorem allows construction of an explicit test function for showing $p$ to be not Bayes. I suspect that the weaker condition
\newline
$\lim_{R \to \infty}\{ \sum_{ij} s_is_j[\int_1^R \textstyle { 1 \over p}V^{-1} r^{1-d}dr]^{-1}ds =0$ is also necessary.
It may be that the condition is also sufficient. 
\noindent A similar condition for the recurrence of diffusion processes is given in [Ic78]. 
\newline\newline
Brown[Br71] studies the admissibility of estimates for the normal location 
problem in $d$ dimensions in which it is assumed that the data $x$ are gaussian with unknown mean and identity covariance matrix. 
He shows that an estimate corresponding to the marginal 
density of the data $ p(x)$  is admissible if $ p_i =\partial _i \log p=\textstyle{\partial 
\over {\partial x_i }}\log p$ is bounded and if
\begin{equation}
  \int_1^\infty {[\int_{s\in 
S} p} ds]^{-1}r^{1-d}dr=\infty .  
\end{equation}
Brown, Theorem 6.4.4,  also shows that an estimate corresponding to the marginal density $ p(x)$  
is admissible only if 
\begin{equation}
 \int_{s\in S} [ \int_1^\infty {p^{-1}} r^{1-d}dr]^{-1}ds=0. 
\end{equation}
\noindent The asymptotic version requires data $X_n \sim N(\theta,I/n)$, 
with $n \to \infty$. 
Theorem 4 implies (40): a prior $ p$  is Bayes  only if 
\begin{equation}
 \int_{s\in S} [ \int_1^\infty {p^{-1}}r^{1-d}dr]^{-1}ds=0 
\end{equation}
 or equivalently, almost everywhere on $S$
\begin{equation} 
\int_1^\infty {p^{-1}} r^{1-d}dr=\infty. 
\end{equation}
If the prior density $ p $ is expressed as a density \textit{$ \rho $ } on the polar co-ordinates 
$x= rs, $ the condition simplifies to $ \int_1^\infty {\rho ^{-1}} (r,s)dr=\infty 
$ almost everywhere on $S$. See Strawderman and Cohen[SC71], theorem 4.4.1.
For example, the prior corresponding to $ r^\alpha $ being uniformly 
distributed is Bayes in every dimension for $ \alpha \le 2$  but not Bayes for 
$ \alpha >2.$ 
\newline\newline
 
\noindent\textbf{7 When is V Bayes on bounded D?}
\newline \newline
Let D be a bounded domain with boundary in $C^{2,\alpha }$ .  Let $ \nu (s),s\in 
\partial D$  denote the outward pointing normal at a point $s$ on the boundary 
of $D$, assumed defined almost everywhere in $ ds$, lebesgue measure on the boundary. It will be assumed 
that, for almost all $s \in \partial D$, the inward pointing normal$ \{s-u\nu (s),0<u<\varepsilon\}$ lies in 
$D$ for $\varepsilon$ small enough. 

\textbf {Theorem 5}. The covariance matrix $pV$ is Bayes  only if 
\begin{equation}
\mathop {\lim }\limits_{\varepsilon \to 0} \int_{s\in \partial D} \nu_i \nu _j 
[\int_0^\varepsilon  {1 \over p} V_{ij}^{-1} (s-u\nu )du]^{-1}ds=0.
\end{equation}

\noindent This is proved similarly to theorem 4, using test functions that are constant on the inward pointing normal segments.
\newline

\noindent It may happen that, for each $s$, the normal vector $ \nu 
(s)$ at $s$ is the limit of some eigenvector of $pV$ as $x \in D \to $s$ \in \partial D$, (the \textit{normal eigenvector} case) in which case the 
condition in the theorem simplifies to
\begin{equation}
  \int_0^\varepsilon {\nu _i \nu _j \over p} V_{ij}^{-1} 
(s-u\nu )du=\infty \mbox{ almost all } s. 
\end{equation}
We will say that the integral condition 
{\textit {fails} at s if $ \int_0^\varepsilon {\nu _i \nu _j \over p} V_{ij}^{-1} (s-u\nu )du<\infty .$ 
The theorem now states that $pV$ is admissible only if the integral condition fails on a set of measure $0$.
\newpage

\noindent The admissible $pV$ are those where $pV \to 0$ fast enough near the boundary. If $pV$ is inadmissible,
 we can render it admissible by attenuating $p$ near the boundary. In the \textit{ normal eigenvector } case, 
let $D_\varepsilon$ consist of those points $x \in D$ within $\varepsilon$ of the boundary, and suppose that each such point 
is closest to a unique boundary point $s(x)$. Each such point may be written $x=s(x)-u(x)\nu(x)$ for some $u$. 

\noindent Let $a$ be an \textit{ attenuation} factor defined at each point in $D$ by:
\begin{eqnarray}
a(x) &=& 1 \mbox{ for } x \in D-D_\varepsilon, \nonumber \\
a(x) &=& 1 \mbox{ for } x \in D_\varepsilon, \mbox{ integral condition holds at } x(s),\nonumber \\
\mathop {g(x )} &=& {\partial \over {\partial u}}(\int_0^{u(x)}{[\nu _i \nu _j }{1 \over p} V_{ij}^{-1} 
(s(x)-w\nu(x) )]^{1/2}dw)^2,\nonumber \\
a(x)&=& min[1,1-(1-{g(x) \over g(\varepsilon})^3 ]\mbox{ for } x \in D_\varepsilon, \mbox{ integral condition fails at } x(s). \nonumber
\end{eqnarray}

\noindent The proposed attenuation factor  will be 1 except near boundary points $ s$  where the integral 
condition fails, where it will approach zero. With the prior $ap$, the integral condition becomes 
\begin{equation}
\int_0^1 {\nu _i \nu _j }{1 \over ap} V_{ij}^{-1} (s-u\nu )du \ge {1 \over 6} g(\varepsilon)\int_0^\varepsilon 
{\partial \over {\partial u}}[log(\int_0^{u}{[\nu _i \nu _j }{1 \over p} V_{ij}^{-1} 
(s-w\nu) dw]du=\infty 
\end{equation}
 \newline \newline

\noindent\textbf{9 One dimension}
\newline \newline
For the one dimensional parameter $x$ on $ D=(a,b)$ with variance $V$, Brown's 
condition implies that $ pV$  is admissible if and only if
\begin{equation}
\int_{a} {(V p)}^{-1}=\int^b {(V p)}^{-1}=\infty .
\end{equation}

\noindent Since a smooth monotone transformation renders $pV$ equal to $1$ on $D=(a,b)$, 
an equivalent result is that there exists a non-zero differentiable test function $w$ on $D$
such that $\int_a^b(w'+{1 \over 2} w^2 )\le 0$ if and only if either $a$ or $b$ are finite.
\newline

\noindent Jeffreys' density $J=V^{-{1 \over 2}}$ is admissible on $D=(a,b)$ if and only if
\begin{equation}
\int_{a} J=\int^b J=\infty,
\end{equation}

\noindent which means that Jeffreys must be "improper" in both tails to be admissible.
I take a certain delight in this impropriety, because although "improper" priors abound in decision theory and in Bayesian analysis, 
they remain objects of suspicion.  See for example the excellent review in [KW96]. However, in decision theory, the prior appears only
when multiplied with a loss function, which may be arbitrarily scaled, so improper priors form a natural part of the range of procedures we need to study. Asymptotically, the prior appears only as a product with the covariance matrix in admissibility questions, and again it makes no sense to constrain priors to be improper. In the Jeffreys' case, the admissibility of the product $pV$ requires that Jeffreys  be improper in the tails.
\newpage
\noindent The pearson correlation coefficient computed for n bivariate 
normal observations with true correlation $ \rho $  has asymptotic variance 
$ 1/(1-\rho ^2)^2$ . 

\noindent Thus a prior $ p$ on $ D=(-1,1)$  is admissible if and only if 
\begin{equation}
 \int_{-1} {((1-\rho ^2} )^2p)^{-1}=\int^1 {((1-\rho ^2} )^2p)^{-1}=\infty .
\end{equation}
 For priors of form $ p=1/(1-\rho ^2)^\alpha $, the prior p is admissible if 
and only if $ \alpha \le 1$ . Thus if we wish to skirt the edge of 
inadmissibility, we might use $ p=1/(1-\rho ^2)$ .
\newline \newline

\noindent\textbf{10 Invariant admissible priors}
\newline\newline

Since the Kullback-Leibler loss function does not change under smooth transformation of the parameter space, 
differences between the asymptotic risk functions for two priors are also invariant under such transformations.
We are free to transform to a convenient $p,V,D$ in deciding admissibility problems.
If a transformation $T$ takes $p,V,D$ into say $T(p),T(V),T(D)$, then the admissibility of $PV$ in $D$ equals the admissibility
of $T(p)V(p)$ in $T(D)$.  A prior $p$ is \textit{ relatively invariant} if $p(Tx)J=cp(x),x \in D$, where $J$ is the Jacobian of the
the transformation $x \to Tx$ and when $T$ is one to one $D \to D$ such that $TV(Tx)=CV(x)$.
\newline\newline
\noindent For example if $D=R^d-{0},V=I$, arbitrary rotations and scalings leave $D$ invariant, and change the covariance by a constant, so
the only invariant priors are of form $p=r^\alpha, r^2=\sum_i x_i^2$. From Brown's condition, the prior p is admissible if
 \begin{equation}
\int_0 r^{1-d-\alpha}dr =\int^\infty r^{1-d-\alpha}=\infty ,
\end{equation}
which occurs only when $\alpha=d-2$. In this case there is a single admissible invariant prior $p=r^{d-2}$. This prior,
discussed in [Br71] and [SC71], corresponds to $r^2$ being uniform over $D$.
\newline\newline
\noindent If $D =\{x|R_1 < r=|x| < R_2\}$, the invariant transformations are rotations, which require that an invariant $p$ depends only 
on $r$. Admissibility requires 
\begin{equation}
\int_{R_1}{1 \over p}dr=\int^{R_2} {1 \over p}dr=\infty.
\end{equation}
Admissibility is achieved by $p(x)=\min_{y \in \partial D}|x-y|$.
\newline\newline
Although invariance considerations no longer always apply, the above solution can be extended to general
bounded $D$ with $V=I$, namely $p(x)=\min_{y \in \partial D}|x-y|$. For general $D,V$, define $|x-y|$ as the path
length between points $x,y \in \bar D$ in the metric $d(x,y)=(x-y)'V^{-1}(x-y)$. Then $D$ is bounded in this metric if all 
paths have finite length, and we again define $p(x)=\min_{y \in \partial D}|x-y|$.
We offer this merely as a suggestion for an admissible prior that flirts with inadmissiblity near the boundaries, and is consistent
under transformations of the data.
\newline\newline
\newpage
\noindent\textbf{11 A mixture model}
\newline \newline
Suppose that $ y_n $  is a sample of size $ n$  from the normal 
mixture
\begin{equation}
 Y=Z+(1-B(q))x_1 +B(q)x_2 
\end{equation}
  where $ Z\sim N(0,1),B(q)\sim \mbox{Bernoulli with mean } q, x_1 >0, x_2 >0, q=\frac{x_1 }{x_1 +x_2 }.$ 

\noindent The parameter $ x=(x_1 ,x_2 )$  lies in the domain $ D=\{x_1 > 0,x_2 > 0\}. 
$ 

\noindent The density of a single observation $ y$  is 
\begin{equation}
f(y)=\{x_2 \phi (y+x_1 )+x_1 \phi (y-x_2 )\}/(x_1 +x_2 ).
\end{equation}

\noindent The asymptotic variance V is the inverse of the information matrix of expected values of 
products of the score functions:
\begin{equation}
l_1 =-\frac{1}{x_1 +x_2 }+\frac{\phi (y-x_2 )}{(x_1 +x_2 )f}-\frac{(y+x_1 
)x_2 \phi (y+x_1 )}{(x_1 +x_2 )f}
\end{equation}
\begin{equation}
l_2 =-\frac{1}{x_1 +x_2 }+\frac{\phi (y+x_1 )}{(x_1 +x_2 )f}+\frac{(y-x_2 
)x_1 \phi (y-x_2 )}{(x_1 +x_2 )f}
\end{equation}
\begin{equation}
L_{ij} =\int l_i l_j fdy
\end{equation}
\begin{equation}
V=L^{-1}
\end{equation}

\noindent Asymptotic admissibility for the prior $p$  is determined by behaviour of $\textit{Vp}$ near 
the boundaries. Let $x_1 =rs_1 ,x_2 =rs_2,s_1^2 +s_2^2 =1.$ 
\begin{eqnarray}
  x_1 \to 0: L_{11} &\to& (\exp (x_2^2 )-1-x_2^2 )/x_2^2, L_{12} \to 0, L_{22} \to 0, 
\label{line1}\\
 x_2 \to 0: L_{22} &\to& (\exp (x_1^2 )-1-x_1^2 )/x_1^2 , L_{12} \to 0, L_{11} \to 0,
\label{line2}\\ 
 r\to \infty : L_{11} &\to& s_2 /(s_1 +s_2 ),L_{12} \to 0, L_{22} \to s_1 /(s_1 +s_2 ). 
\label{line3}
\end{eqnarray}

\noindent At the boundary $ x_1 =0$ the normal is an eigenvector at all points $(0,x_2)$, and the integral 
condition for admissibility for that boundary is $ \int_0^1 {1 \over p}L_{11} dx_1 =\infty $ 
almost all $ x_2 $ which reduces to $ \int_0^1 {1 \over p} dx_1=\infty$ almost all $ x_2 $.
Similarly, the condition for admissibility on the boundary $ x_2 =0$ is 
$ \int_0^1 {1 \over p} dx_2=\infty$ almost all $ x_1 $.

\noindent For the infinite ``boundary'' $r \to \infty$, the integral condition for admissibility is 
\begin{equation}
  \lim _{R\to \infty } \int_{s\in S} s_i s_j [\int_1^R {1 \over p} V_{ij}^{-1} (rs)r^{-1}dr]^{-1}ds=0 ,
\end{equation} 

 \noindent where S is the intersection of the 
boundary of the unit circle and the upper right quadrant. Using the behavior 
of $L$ as $r \to \infty$, this condition becomes 
\begin{equation}
\int_1^R {1 \over rp(sr)} dr \to \infty \mbox{ almost all } s \in S.
\end{equation}

\noindent Choosing a prior p to make pV admissible requires that 
\begin{equation}
\int_0^1 {1 \over p} dx_1 =\int_0^1 {1\over p}dx_2  ,\int_1^\infty {1 \over pr}dr=\infty .
\end{equation}
Roughly, we need that $p$ be of order $ x_1 $  near $ x_1 =0$ , 
of order $ x_2 $  near $ x_2 =0$ ,and of order $log(r)$ near  $ r=\infty $ .
For example, $ p=\frac{x_1 x_2 }{(x_1 +x_2 )^2}$  will do the job, as will many 
other priors with the correct behavior near the boundary. The uniform is inadmissible because it fails at 
$x_1=0$ and $x_2=0$.
\newpage

\begin{figure}[t] 
  \centering
  \includegraphics[bb=56 134 556 658,width=5.67in,height=5.94in,keepaspectratio]{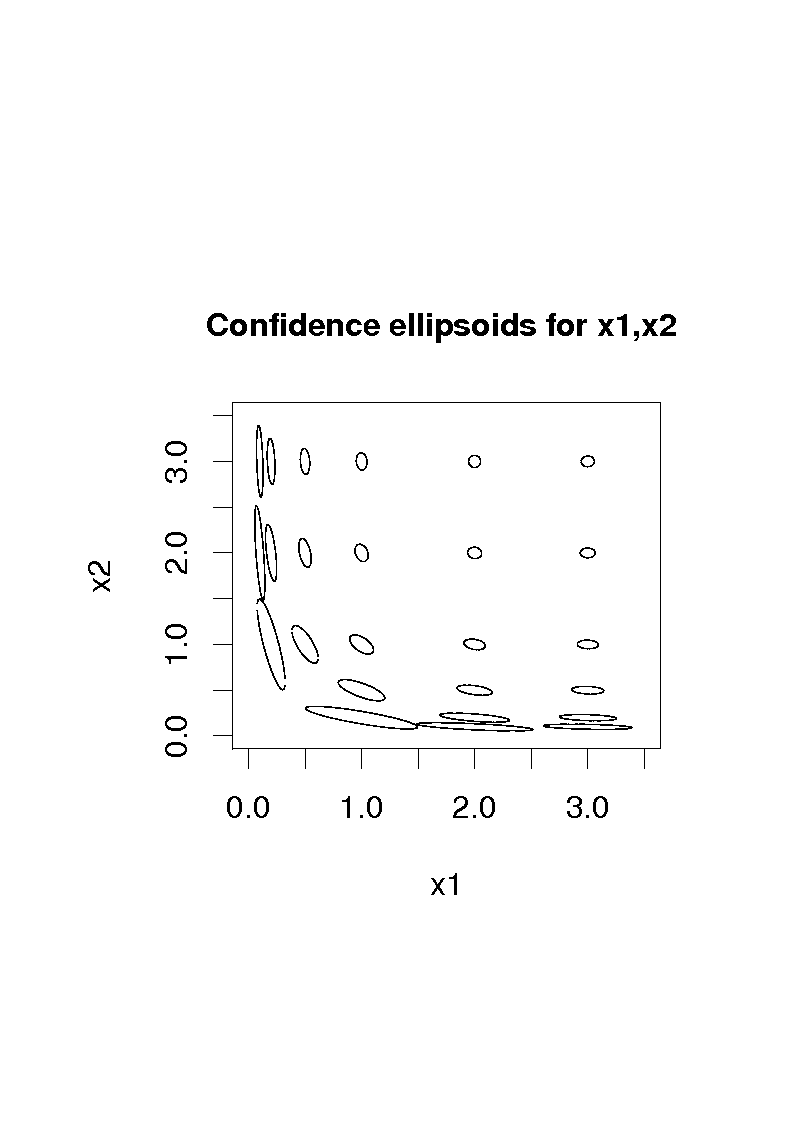}
  \caption{Confidence Ellipses}
  \label{fig:confidence}
\end{figure}

The plot of confidence ellipses when 1000 points are sampled from the 
mixture model shows how the boundaries affect asymptotic admissibility. For 
the boundaries $x1=0$ and $x2=0$, the asymptotic variances orthogonal to the boundary in 
fact approach a positive limit; thus the integral of the inverse variances 
up to the boundary is positive rather than infinite, and the uniform density 
is therefore inadmissible. For the boundary at infinity, the variances are 
bounded away from infinity, so the integral of the inverse variance is 
infinite, and this boundary is admissible for a uniform prior.

 \newpage
\noindent \textbf{12 A prior beating the uniform}
\newline \newline
It is of interest to exhibit a prior with asymptotic risk everywhere smaller 
than an inadmissible prior such as the the uniform in this problem.  Brown's condition exhibits
a prior satisfying $ \textstyle{\partial \over {\partial 
x'}}[V\textstyle{{\partial p} \over {\partial x}})=0. $ The asymptotic risk of 
$p$, relative to the uniform, is $ -\textstyle{1 \over 2}(\textstyle{{\partial 
p} \over {\partial x}})'V\textstyle{{\partial p} \over {\partial x}}/p$ . 
There are many solutions to the elliptic differential equation, depending on 
boundary values of $p$. The solutions are not necessarily admissible. 

We have computed solutions for the discrete approximation where $ x_1 ,x_2 $  
each lie in the grid 0.1, 0.2,{\ldots}10. The boundary values for $p$ are 
$ p=\min (\frac{4x_1 x_2 }{(x_1 +x_2 )^2},x_1 x_2 ). $  We set these values so 
that p will satisfy the conditions for admissibility at the different 
boundaries. The following prior is obtained by using a relaxation method to 
solve the finite difference form of the differential equation; at the 
solution, the finite difference expressions are everywhere less than .01. 
A similar prior was developed in [Em02].
\begin{figure}[p] 
  \centering
  \includegraphics[bb=0 0 840 525,width=5in,height=4in,keepaspectratio]{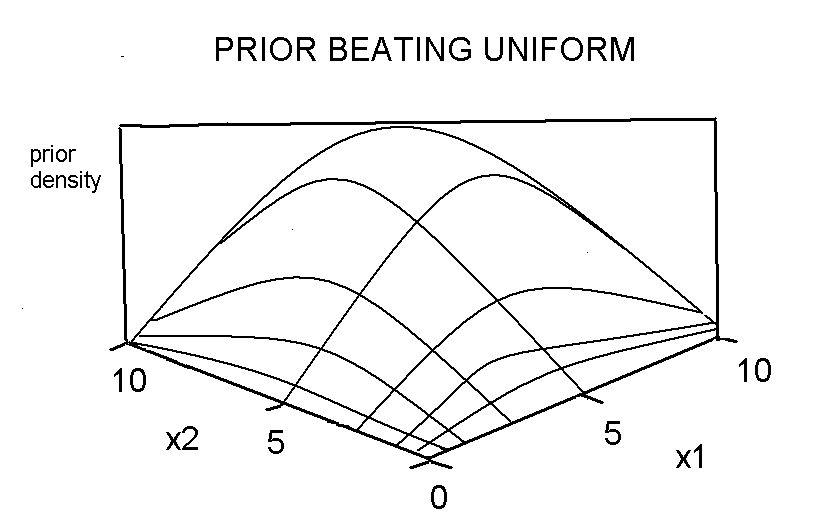}
  \caption{Prior beating uniform}
  \label{fig:priorpaint}
\end{figure}

It will be noted that the prior density approaches zero at the lower and 
left boundary, but not at the other two boundaries, as required by the 
admissibility conditions.

\begin{figure}[p] 
  \centering
  \includegraphics[bb=0 0 840 525,width=5.67in,height=3.54in,keepaspectratio]{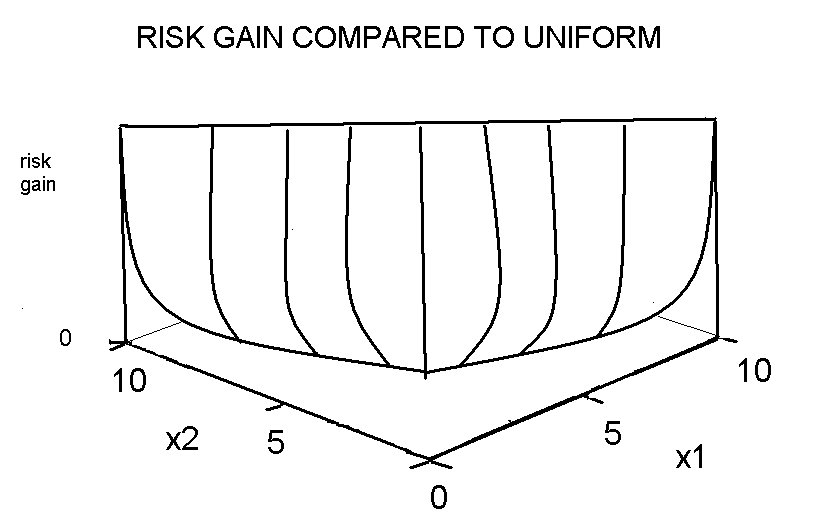}
  \caption{Risk Gains compared to Uniform}
  \label{fig:riskgain}
\end{figure}

The risk gains against the uniform are everywhere positive ( as required by 
the theory), but are far greater near the low $x1$ and $x2$ boundaries. This is 
to be expected, because the prior is made admissible by changes near the boundaries, so
that larger improvements in the risk should occur there. 
\newline \newline
\newpage 
\noindent \textbf{13 Acknowledgement}

Although Bill Strawderman has not seen this particular manuscript, he has over the years seen many an
unsolicited manuscript from me on similar topics, and has always responded generously with his time, his experience, and his
vast penetrating knowledge of this subject. Thanks Bill for the many kind assistances, for the numerous balloon-bursting 
counterexamples, and for the splendid life example of great imagination, clear thinking, clear exposition, and genial princely decency.


\end{document}